%&amstex          
\input amstex\documentstyle{amsppt}  
\pagewidth{12.5cm}\pageheight{19cm}\magnification\magstep1
\topmatter
\title{Rigid strata in a reductive group}\endtitle
\author G. Lusztig\endauthor
\address{Department of Mathematics, M.I.T., Cambridge, MA 02139}\endaddress
\thanks{Supported by NSF grant DMS-2153741}\endthanks
\endtopmatter   
\document

\define\mpb{\medpagebreak}

\define\hW{\hat W}

\define\frl{\forall}

\define\si{\sim}

\define\sqc{\sqcup}

\define\qua{\quad}

\define\op{\oplus}
   
\define\part{\partial}

\define\n{\notin}

\define\m{\mapsto}
\define\do{\dots}

\define\sub{\subset}    

\define\T{\times}

\define\nl{\newline}
\redefine\i{^{-1}}

\define\un{\underline}

\define\Ad{\text{\rm Ad}}
\define\Hom{\text{\rm Hom}}

\define\ind{\text{\rm ind}}

\redefine\c{\chi}

\define\io{\iota}

\define\r{\rho}
\define\s{\sigma}

\redefine\l{\lambda}

\define\Si{\Sigma}

\define\kk{\bold k}

\define\QQ{\bold Q}

\define\cl{\Cal L}

\define\co{\Cal O}
\define\cp{\Cal P}

\define\ct{\Cal T}
\define\cu{\Cal U}

\head Introduction\endhead
\subhead 0.1\endsubhead
Let $G$ be a connected reductive group over an algebraically
closed field $\kk$. Let $\l(G)$ be the set of subgroups
$L\sub G$ which are Levi subgroups of parabolic
 subgroups of $G$.
 Let $\cu(G)$ be the set of unipotent
conjugacy classes in $G$. We say that $C\in\cu(G)$ is rigid
(see \cite{S82, 7.9})
if $C$ cannot be obtained by the induction procedure of
\cite{LS79} from a unipotent class of $L$ with $L\in\l(G)$,
$L\ne G$. Let $\cu_{rig}(G)$ be the subset of $\cu(G)$
consisting of rigid unipotent classes.

In \cite{L15} we have defined a partition of $G$ into finitely
many strata; each stratum is locally closed \cite{Ca20} and
a union of $G$-conjugacy classes
of fixed dimension. The set $Str(G)$ of strata of $G$ can be
viewed as an enlargement of $\cu(G)$ (a unipotent class of $G$
is contained in exactly one stratum). In this paper we extend
the notion of rigid unipotent class and introduce a subset
$Str_{rig}(G)$ of $Str(G)$ whose elements are called
{\it rigid } strata.
Namely, we say that a stratum $\Si$ of $G$ is rigid if it
is a finite union of orbits of $G\T Z_G^0$ on $G$, where the
first factor acts by conjugation and the second factor (the
connected centre of $G$) acts by translation.
This definition is quite different from the definition of rigid
unipotent classes. To remedy this, we shall also give a
second definition of rigidity of a stratum based on a notion of
induction from $L$ to $G$ of a stratum of $L$ where
$L\in\l(G)$; we show that it is equivalent to the first.
Using this second definition we can relate the notion
of rigid stratum to that of rigid unipotent classes in various
versions of $G$ in different characteristics. We can use this
relation to give a classification of rigid strata of $G$
based on the classification of rigid unipotent classes in
\cite{S82}; we will do this explicitly for $G$ of type $E_6,E_7,E_8$
(see 2.5, 2.6).

\subhead 0.2\endsubhead
Let $W=W_G$ be the Weyl group of $G$. 
Let $cl(W)$ be the set of conjugacy classes in $W$. In \cite{L15}
a surjective map $cl(W)@>>>Str(G)$ is defined and in \cite{L22}
it is shown that each fibre of this map contains a unique element
which as non-elliptic as possible; assigning to a stratum this
unique element defines an injective map $Str(G)@>>>cl(W)$ whose
image is denoted by $CL(W)$. Thus we obtain a bijection $Str(G)@>\si>>CL(W)$.
The image of $Str_{rig}(G)$ under this bijection is denoted by $CL_{rig}(W)$; we say
that this is the set of rigid conjugacy classes in $W$.

\subhead 0.3\endsubhead
{\it Notation.} 
 Let $\hW=\hW_G$ be the set of isomorphism classes of
irreducible $\QQ[W]$-modules. Let $\un\hW_G$ be the image of
the imbedding $\io^G:\cu(G)@>>>\hW$ given by the Springer
correspondence (extended in \cite{L84} to small characteristic).
For $x\in G$ we write $x=x_sx_u=x_ux_s$ where $x_s$ (resp.
$x_u$) is semisimple (resp. unipotent).
For $s\in G$ semisimple let $Z^0_G(s)$ be the connected
centralizer of $s$ in $G$ and let $W_G(s)$ be its Weyl group.

\head 1. Rigid strata\endhead
\proclaim{Lemma 1.1} Let $V$ be a finite dimensional $\kk$-vector space.
Let $A:T@>>>GL(V)$ be a rational representation of a torus $T$ over
$\kk$. Let $\s\in GL(V)$ be semisimple and such that $A(t)\s=\s A(t)$
for all $t\in T$. There exists an open dense subset $\co$ of $T$ such
that the following holds. For any $z\in\co$ and any subspace $V'$ of
$V$ such that $\s A(z)(V')=V'$, we have $\s(V')=V'$ and $A(t)(V')=V'$
for all $t\in T$.
\endproclaim
For $c\in\kk^*$ let $V_c=\{v\in V;\s(v)=cv\}$.
For $\c\in X=\Hom(T,\kk^*)$ let 
$V^\c=\{v\in V;A(t)v=\c(t)v \qua\frl t\in T\}$.
Then $V=\op_{(c,\c)\in J\T J'}V_c\cap V^\c$ where
$J\sub\kk^*,J'\sub X$ are finite.

Let $\co$ be the set of all $t\in T$ such that
$$c\c(t)\ne c'\c'(t)\text{ for any }(c,\c)\ne(c',\c')\text{ in }J\T J'$$
that is, $(\c\c'{}\i)(t)\ne c'c\i$ for any $(c,\c)\ne(c',\c')$ in
$J\T J'$. If $\c=\c'$ this condition is clearly
satisfied. If $\c\ne\c'$, this means
$(\c\c'{}\i)(t)\ne c'c\i$ for $c,c'$ in $J$ that is,
$$t\n\cup_{\c\ne\c' \text{ in } J',c,c'\text{ in } J}(\c\c'{}\i)\i(c'c\i);
$$

thus, $\co$ is the complement of the union of finitely many
codimension $1$ subvarieties of $T$. In particular $\co$ is open
dense in $T$.

Assume now that $z\in\co,V'$ are as in the lemma.
Since $\s A(z)(V')=V'$ we have $V'=\op_{a\in\kk^*}V'_a$ where
$V'_a=\{v\in V';\s A(z)v=av\}$.
Let $v\in V'_a$. We have $v=\sum_{c,\c}v'_{c,\c}$ where
$v'_{c,\c}\in V_c\cap V^\c$ satisfy $c\c(z)=a$.
Applying $\s A(z)$, we get
$$a\sum_{c,\c}v'_{c,\c}=\sum_{c,\c}c\c(z)v'_{c,\c}$$
that is,
$$\sum_{c,\c}(a-c\c(z))v'_{c,\c}=0,$$
so that $(a-c\c(z))v'_{c,\c}=0$ for any $(c,\c)\in J\T J'$.
Hence for any $c,\c$ such that $v'_{c,\c}\ne0$, we have $c\c(z)=a$.
Thus, $v=\sum_{c,\c;c\c(z)=a}v'_{c,\c}$.
Since $z\in\co$, this sum has at most one nonzero term, so that
$v=0$ or $v\in V_c\cap V^\c$ for a unique $c,\c$. In particular,
for some $c,\c$ we have $\s v=cv\in V'$ and $A(t)(v)=\c(t)v\in V'$
for all $t\in T$, so that $\s(V')\sub V'$ and $A(t)(V')\sub V'$ for
all $t\in T$. The lemma is proved.

\subhead 1.2\endsubhead
Let $L\in\l(G)$, $x\in L$. Let $V=Lie(G)$, $T=Z_L^0$. Define
$A:T@>>>GL(V)$ by $A(t)=Ad(t)$. Define $\s\in GL(V)$ by $\s=\Ad(x_s)$.
Let $\co=\co_x$ be the open dense subset of $T$ defined in Lemma 1.1
in terms of $V,T,A,\s$.

\proclaim{Lemma 1.3} We preserve the setup of 1.2. We set $x_s=s$.
For any $z\in\co_x$ we have $Z^0_G(x_sz)=Z^0_L(s)$.
\endproclaim
Clearly, we have $Z^0_L(s)\sub Z^0_G(sz)$ for any $z\in T$.
It remains to show that if $z\in\co_x$, then
$Lie(Z^0_G(sz))\sub Lie(Z^0_L(s))$. Let $\cl$ be a line in
$Lie(Z^0_G(sz))$. We must show that $\cl\sub Lie(Z^0_L(s))$.
We have $\s A(z)(\cl)=\cl$. By Lemma 1.1 
we then have $\s(\cl)=\cl$ and $A(t)(\cl)=\cl$ for all $t\in T$. The
last equality shows that $\cl\sub Lie(L)$ and the equality
$\s(\cl)=\cl$ shows that $\cl\sub Lie(Z^0_G(s)$. It follows that
$\cl\sub Lie(Z^0_G(s))\cap Lie(L)=Lie(Z^0_L(s))$. The lemma is proved.

\subhead 1.4 \endsubhead
Let $x\in G$. We set $s=x_s,u=x_u$. 
Let $\r^G_u=\io^G(C)\in\hW^G$ (see 0.3) where $C$ is the
$G$-conjugacy class of $u$. Following \cite{L15, 2.2} we set
$$\r^G_x=j_{W_G(s)}^W(\r^{Z^0_G(s)}_u)\in\hW.\tag a$$
Here $j_{W_G(s)}^W$ is the truncated induction defined in
\cite{LS79}; $W_G(s)$ is regarded as a reflection subgroup of
$W$ as in \cite{L15, 2.1}.

We now assume that $L\in\l(G)$, $x\in L$ (so that $s\in L,u\in L$);
$T=Z_L^0$, $\co_x\sub T$ are as in 1.2.
For any $z\in\co$, $\r^G_{xz}\in\hW$, $\r^L_x\in\hW_L$ are defined
as above. We show:
$$\r^G_{xz}=j_{W_L}^W(\r^L_x)\in\hW.\tag b$$
An equivalent statement is
$$j_{W_G(sz)}^W(\r^{Z^0_G(sz)}_u)=
j_{W_L}^W(j_{W_L(s)}^{W_L}(\r^{Z^0_L(s)}_u).$$
Using the transitivity of truncated induction, this is the same as
$$j_{W_G(sz)}^W(\r^{Z^0_G(sz)}_u)=(j_{W_L(s)}^W(\r^{Z^0_L(s)}_u).
\tag c
$$

By Lemma 1.3 we have $Z^0_G(sz)=Z^0_L(s)$. It follows that
$W_G(sz)=W_L(s)$ so that (c) holds. We see that (b) holds.

\subhead 1.5\endsubhead
Recall from \cite{L15} the definition of strata of $G$. Two elements
$g,g'$ of $G$ are said to be in the same stratum of $G$ if
$\r^G_g=\r^G_{g'}$ as objects of $\hW$. Let $Str(G)$ be the set of
strata of $G$. Note that $\Si\m\r^G_\Si=\r^G_g$ for any $g\in\Si$ is a
well defined injective map $Str(G)@>>>\hW$;
let $\hW_G^*$ be the image of this map. Thus, $\Si\m\r^G_\Si$ is a
bijection
$$Str(G)@>\si>>\hW_G^*;\tag a$$
in particular $Str(G)$ is finite.

\subhead 1.6\endsubhead
Now let $L\in\l(G)$. For any $x\in L$ let $\co_x\sub Z_L^0$ be as in
1.3.

(a) {\it Let $\Si'$ be a stratum of $L$. There is a unique stratum $\Si$
of $G$ such that for any $x\in\Si'$ and any $z\in\co_x$ we have
$xz\in\Si$.}
\nl
The uniqueness of $\Si$ is obvious. To prove its existence it is enough
to show that if $x,x'$ in $L$ satisfy $\r^L_x=\r^L_{x'}$ and if
$z\in\co_x,z'\in\co_{x'}$ then $\r^G_{xz}=\r^G_{x'z'}$. This
follows immediately from 1.4(b).

We say that $\Si$ is obtained from $\Si'$ by induction from $L$.
 From 1.4(b) we see that

(b) $\r^G_\Si=j_{W_L}^W(\r^L_{\Si'})$.

\subhead 1.7\endsubhead
Recall from \cite{L84, 3.1} the definition of pieces of $G$.
Let $L\in\l(G)$ and let $S$ be the inverse image under $L@>>>L/Z_L^0$
of an isolated conjugacy class (see \cite{L84, 2.6}) of $L/Z_L^0$.
For such $(L,S)$ we define $[L,S]$ to be the union of $G$-conjugates of
the subset $\{x\in S;Z^0_G(x_s)\sub L\}$ of $S$. The subset $[L,S]$ is
said to be a piece of
$G$. It depends only on the $G$-conjugacy class of $(L,S)$. The pieces
of $G$ form a partition of $G$ into finitely many subsets. We show:

(a) {\it If $[L,S]$ is a piece then $[L,S]$ is contained in the
stratum of $G$ obtained by inducing from $L$ to $G$ the stratum of
$L$ containing $S$.}
\nl
Let $x\in S$ be such that  $Z^0_G(x_s)\sub L$.
It is enough to show that for any $z\in\co_x$ (notation of 1.2)
we have $Z^0_G((xz)_s)\sub L$ that is $Z^0_G(x_sz)\sub L$.
By Lemma 1.3 we have $Z^0_G(x_sz)=Z^0_L(x_s)$ and in particular
$Z^0_G(x_sz)\sub L$. This proves (a).

From (a) we see that any stratum of $G$ is a union of pieces of $G$
(this was stated without proof in \cite{L15}).

\subhead 1.8\endsubhead
The $G\T Z_G^0$ action on $G$ (see 0.1) preserves each piece of $G$.
We show:

(a) {\it A piece $[L,S]$ contains only finitely many $G\T Z_G^0$-orbits if
and only if $L=G$ (in which case it is a single $G\T Z_G^0$-orbit,
namely $S$).}
\nl
Replacing $G$ bt $G/Z_G^0$, we can assume that $Z_G^0=\{1\}$.
The ``if'' part is obvious. Conversely, assume that
 $[L,S]$ contains only finitely many $G$-conjugacy classes.
Let $x\in S$. As in the proof of 1.7(a) we have
$xz\in [L,S]$ for any $z\in\co_x$ (notation of 1.2).
Thus $x\co_x$ meets only finitely many $G$-conjugacy classes.
Taking semisimple parts we see that
$x_s\co_x$ meets only  finitely many $G$-conjugacy classes.
Let $\ct$ be a maximal torus of $L$ containing $x_s$; we
have $\co_x\sub Z_L^0$ hence $\co_x\sub\ct$ so that $x_s\co_x\sub\ct$.
Now any $G$-conjugacy class of $G$ intersects $\ct$ in a finite set.
It follows that $x_s\co_x$ is a finite set. Since $\co_x$ is open dense
in $Z_L^0$ it follows that $Z_L^0=\{1\}$ so that $L=G$. This proves (a).

\proclaim{Proposition 1.9} Let $\Si$ be a stratum of $G$. Then $\Si$
is rigid (see 0.1) if and only if it satisfies (i) below.

(i) for any $L\in\l(G)$ with $L\ne G$ and any $\Si'\in Str(L)$, $\Si$
is not the stratum induced from $L,\Si'$.
\endproclaim
As in 1.8, we can assume that $Z_G^0=\{1\}$.
Assume first that $\Si$ does not satisfy (i), but $\Si$ is a finite union of
$G$-conjugacy classes. We can then find $L\in\l(G)$
with $L\ne G$ and $\Si'\in Str(L)$ such that for any $x\in\Si'$ and
any $z\in\co_x$ we have $xz\in\Si$. Then 
$x\co_x$ meets only finitely many $G$-conjugacy classes.
Repeating the argument in the proof of 1.8(a) we see that $L=G$, a
contradiction. We see that if $\Si$ is a finite union of $G$-conjugacy
classes, then $\Si$ satisfies (i).
Next we assume that $\Si$ is a union of infinitely many
$G$-conjugacy classes. Since $\Si$ is a finite union of pieces, there
exists at least one piece, say $[L,S]$, which is contained in $\Si$ and
is a union of infinitely many $G$-conjugacy classes. 
From 1.8(a) we see that $L\ne G$ and from 1.7(a) we have that $\Si$
is induced from a stratum of $L$. Thus $\Si$ does not satisfy (i). This
proves the proposition.

\subhead 1.10\endsubhead
Let $Str_{rig}(G)$ be the subset of $Str(G)$ consisting of
rigid strata.

We say that $E\in\hW_G^*$ is rigid if for any $L\in\l(G)$ with $L\ne G$
and any $E'\in \hW_L^*$, we have $E\ne j_{W_L}^W(E')$.
Let $(\hW_G^*)_{rig}$ be the subset of $\hW_G^*$ consisting of
rigid elements. Using 1.6(b) we see that under the bijection
$Str(G)@>\si>>\hW_G^*$ in 1.5(a), $Str_{rig}(G)$ corresponds to
$(\hW_G^*)_{rig}$.

\head 2. Relation with rigid unipotent classes\endhead
\subhead 2.1\endsubhead
Let $\cp=\{2,3,5,7,\do\},\bar\cp=\cp\sqc\{0\}$. For
$r\in\bar\cp$ let $\kk_r$ be an
algebraically closed field of characteristic $r$ and let $G^r$ be a
connected reductive group over $\kk_r$ of the same type as $G$ (hence
with the same Weyl group $W$); we can assume that for some $r\in\bar\cp$
we have $\kk=\kk_r,G=G^r$. From \cite{L15} it is known that for any
$r'\in\bar\cp$ we have
$$\hW_{G^{r'}}^*=\cup_{r\in\bar\cp}\un\hW_{G^r}\sub\hW;\tag a$$
in particular, $\hW_{G^{r'}}^*$ is independent of $r'$ so that it is
the same as $\hW^*_G$. It also follows that $Str(G^{r'})$ can be
indexed in a manner independent of $r'$.

\subhead 2.2\endsubhead
Let $L\in\l(G)$. This is the Levi subgroup of a parabolic subgroup
$P$ of $G$. For each $r\in\bar\cp$ let $L^r\in\l(G^r)$ be such that
$L^r$ is the Levi subgroup of a parabolic subgroup of $G^r$ of the same
type as $P$; we can assume that $W_{L^r}=W_L$ for any $r\in\bar\cp$.
Let $r\in\bar\cp$.
Let $\ind_{L^r}^{G^r}:\cu(L^r)@>>>\cu(G^r)$ be the induction of unipotent
classes from \cite{LS79}.
Under the identifications $\io^{G^r}:\cu(G^r)@>\si>>\un\hW_{G^r}$,
$\io^{L^r}:\cu(L^r)@>\si>>\un\hW_{L^r}$, the map $\ind_{L^r}^{G^r}$
becomes the truncated induction map
$j_{W_L}^W:\un\hW_{L^r}@>>>\un\hW_{G^r}$ (see \cite{LS79}).

\subhead 2.3\endsubhead
Let $r\in\bar\cp$. Recall from \cite{S82, 7.9} that $C\in\cu(G_r)$ is
said to be rigid if for any $L\in\l(G),L\ne G$ and any
$C'\in\cu(L^r)$ we have $C\ne\ind_{L^r}^{G^r}(C')$.

Let $\cu_{rig}(G^r)$ be the subset of $\cu(G^r)$
consisting of rigid unipotent classes.
Let $(\un\hW_{G^r})_{rig}$ be the set of objects of
$E\in\un\hW_{G^r}$ such that for any $L\in\l(G),L\ne G$ and any
$E'\in\un\hW_{L^r}$ we have $E\ne j_{W_L}^W(E')$.
From the results in 2.2 we see that $\io^{G^r}$ restricts to a
bijection $\cu_{rig}(G^r)@>\si>>(\un\hW_{G^r})_{rig}$.

\subhead 2.4\endsubhead
Let $r'\in\bar\cp$. We show:
$$(\hW_{G^{r'}}^*)_{rig}=\cup_{r\in\cp}(\un\hW_{G^r})_{rig}.\tag a
$$
We can assume that $G$ is simple.
Assume first that $G$ is of type $\ne E_8,G_2$. From \cite{L84}, \cite{S85}
it is known that
for any $L\in\l(G)$ and any
$r\in\bar\cp$ we have $\un\hW_{L^r}\sub\un\hW_{L^2}$
hence

(b) $\hW_{L^{r'}}^*=\un\hW_{L^2}$.
\nl
Then (a) follows immediately from the definition.

Assume next that $G$ is of type $G_2$.  From \cite{S85}
it is known that
for any $L\in\l(G)$ and any
$r\in\bar\cp$ we have $\un\hW_{L^r}\sub\un\hW_{L^3}$
hence

(c) $\hW_{L^{r'}}^*=\un\hW_{L^3}$.
\nl
Then (a) follows immediately from the definition.

Now assume that $G$ is of type $E_8$. It is known \cite{L84},
\cite{S85}, that
for any $L\in\l(G)$ with $L\ne G$ and any $r\in\bar\cp$ we have
$\un\hW_{L^r}\sub\un\hW_{L^2}$  hence
$\hW_{L^{r'}}^*=\un\hW_{L^2}$. It is also known \cite{S85}
that for any $r\in\bar\cp$ we have
$\un\hW_{G^r}\sub\un\hW_{G^2}\sqc\{E\}$
where $E\in\hW$ is the unique object in
$\un\hW_{G^3}$ which is not in $\un\hW_{G^2}$; hence

(d) $\hW_{G^{r'}}^*=\un\hW_{G^2}\sqc\{E\}$.
\nl
It is then enough to observe that for any $L\in\l(G),L\ne G$ and any
$E'\in\un\hW_{L^2}$ we have  $E\ne j_{W_L}^W(E')$; indeed,we have
$j_{W_L}^W(E')\in\un\hW_{G^2}$ and $E\in\un\hW_{G^2}$. This proves (a).

\subhead 2.5\endsubhead
Using 2.4(b),(c),(d) we can calculate explicitly the set $(\hW_G^*)_{rig}$
for any simple $G$. We show the results for $G$ simple of type $E_6,E_7,E_8$.
These results are obtained using \cite{L82}, \cite{S85}.
(The elements of $\hW$ are denoted as in \cite{S85}.)

\mpb

Type $E_6$.

$1_{36},6_{25},15_{16},10_9$.

\mpb

Type $E_7$.

$1_{63},7_{46},27_{37},35_{31},15_{28},189_{22},70_{18},280_{17},84_{15}$.  

\mpb

Type $E_8$.

$1_{120},8_{91},35_{74},84_{64},50_{56},210_{52},560_{47},400_{43},448_{39},$ 

$1344_{38},175_{36},1050_{34},972_{32},1400_{29},840_{26},168_{24},420_{20},$

$1344_{19},2016_{19},840_{14},175_{12}$. 

\subhead 2.6\endsubhead
We now describe the set $CL_{rig}(W)$ for $G$ simple of type
$E_6,E_7,E_8$.
(We use 2.5 and the tables in \cite{L15,\S4}). The notation of conjugacy classes in $W$
is as in \cite{C72}.

\mpb

Type $E_6$.

$A_0,A_1,3A_1,2A_2A_1$.

\mpb

Type $E_7$.

$A_0,A_1,2A_1,(3A_1)'',(4A_1)'',A_2+2A_1, 2A_2+A_1,(A_3+A_1)'',A_3+A_2$.

\mpb

Type $E_8$.

$A_0,A_1,2A_1,3A_1,(4A_1)'',A_2+A_1,A_2+2A_1,A_2+3A_1,2A_2+A_1,$

$A_3+A_1,2A_2+2A_1,(A_3+2A_1)'',A_3+A_2,A_3+A_2+A_1,(2A_3)'',D_4+A_2,$

$A_4+A_3,D_5(a_1)+A_2,(A_5+A_1)'',D_5+A_2,(A_7)''$.

\enddocument